\documentclass[conference,letterpaper]{IEEEtran}

\IEEEoverridecommandlockouts

\usepackage{cite}
\usepackage{amsmath,amssymb,amsfonts}
\usepackage{graphicx}
\usepackage{textcomp}
\usepackage{algorithm}
\usepackage{algorithmic}
\usepackage{comment}
\usepackage{soul}
\usepackage{xcolor}
\usepackage[normalem]{ulem} 

\usepackage{geometry}
\usepackage{afterpage}
\AtBeginDocument{%
  \newgeometry{letterpaper,left=54pt,right=54pt,top=72pt,bottom=54pt}
  \afterpage{\restoregeometry}
}

\def\BibTeX{{\rm B\kern-.05em{\sc i\kern-.025em b}\kern-.08em
    T\kern-.1667em\lower.7ex\hbox{E}\kern-.125emX}}

\begin{document}

\title{Optimal Control for Cancer Chemotherapy Using Hybrid Quantum Particle Swarm Optimization}

\author{\IEEEauthorblockN{Bereket Sitotaw Kidane, Md Samiul Haque Motayed, and Shuo Wang\IEEEauthorrefmark{1}}
    \thanks{This work was supported by SW’s faculty Science and Technology Acquisition and Retention (STARs) Program (project ID: AR91084L-51) and NSF CAREER Award (FAIN: 2238269).}
    \thanks{B. S. Kidane and M. S. H. Motayed are Ph.D. students in the Department of Mechanical and Aerospace Engineering, University of Texas at Arlington, Arlington, TX 76019, USA.}
    \thanks{\IEEEauthorrefmark{1}S. Wang is with the Faculty of the Department of Mechanical and Aerospace Engineering, University of Texas at Arlington, Arlington, TX 76019, USA (corresponding author's e-mail: shuolinda.wang@uta.edu).}
}

\maketitle

\begin{abstract}
Optimal control in cancer chemotherapy is challenged by tumor heterogeneity and mutations, which complicate the effectiveness of treatment. Traditional methods, such as Pontryagin's maximum principle (PMP), are often hindered by their reliance on an initial guess for the costate equation, affecting accuracy and convergence. To address these limitations, this work introduces a hybrid Quantum Particle Swarm Optimization (QPSO) method based on regularization. QPSO is employed for global exploration to approximate the optimal control trajectory, followed by a regularization-based refinement to ensure smoothness and consistency with optimality conditions. The Hamiltonian function is used for first- and second-order optimality checks, verifying solution quality. Numerical case studies explore various drug effectiveness functions, demonstrating the role of periodic and localized drug delivery in achieving robust tumor suppression and providing insights into the impact of different drug combinations on optimal chemotherapy strategies.
\end{abstract}


\section{Introduction}

Cancer chemotherapy faces significant challenges due to tumor heterogeneity and mutation dynamics, which often lead to resistance to treatment and reduced therapeutic efficacy. Tumor heterogeneity, characterized by diverse phenotypic and genetic variations among tumor cells, complicates the development of therapies capable of targeting all malignant cells effectively \cite{dagogo2018tumour, diaz2012tumor, lavi2013role, ledzewicz2014optimal}. This diversity enables subpopulations of cancer cells to evade treatment, leading to disease progression. Mutation dynamics further exacerbate this issue by allowing tumor cells to evolve in response to therapeutic interventions, thereby diminishing drug effectiveness over time and increasing the likelihood of treatment failure \cite{gatenby2009adaptive, fu2025intratumoral, greene2014impact}. As highlighted in recent studies, this evolving nature of tumors necessitates dynamic treatment strategies that can respond to tumor evolution in real-time, ensuring sustained therapeutic success \cite{wang2016optimal, hahnfeldt2003minimizing}.

Despite advances in treatment strategies, conventional chemotherapy regimens often fail to account for these dynamic changes, typically relying on fixed-dose protocols or heuristic scheduling methods that do not adapt to the changing characteristics of the tumor. This limitation underscores the need for innovative approaches that integrate mathematical modeling, optimization, and control theory to improve the efficacy of chemotherapy. As shown in the literature, the integration of optimal control theory, such as bang-bang and singular controls, holds the potential to develop more effective treatment protocols by balancing tumor reduction with minimizing side effects \cite{glick2017optimal, ledzewicz2016drug, schattler2015optimal}. Recent advances in computational oncology use mathematical modeling, optimal control, and artificial intelligence for personalized and adaptive chemotherapy, but high‑dimensional, nonlinear tumor dynamics make optimization challenging~\cite{wang2017fixed, wang2018free, rabiei2024long, khudhair2025nonlinear}. These challenges highlight the need for optimization frameworks to navigate complex solution spaces while ensuring smooth and clinically viable drug administration policies.

A major hurdle in chemotherapy optimization is that traditional methods like the Pontryagin's Maximum Principle (PMP) and shooting method struggle with convergence in high‐dimensional tumor models and are highly sensitive to initial conditions\cite{pahnehkolaei2023optimal}. These methods often fail for nonlinear, multimodal tumor dynamics, where small changes in initial guesses yield vastly different outcomes. Robust, adaptive optimization is needed to balance global search and local refinement for feasible, smooth chemotherapy protocols.

To address these challenges, we use a hybrid quantum optimization and regularization framework that integrates QPSO with a regularization-based local refinement using fmincon. Quantum particle swarm optimization (QPSO) builds on classical Particle Swarm Optimization (PSO) by incorporating quantum mechanics-inspired principles, allowing particles to explore the solution space more effectively \cite{yang2004quantum}. QPSO leverages quantum potential fields for efficient global search with less reliance on precise initial conditions—ideal for optimizing chemotherapy in complex, multimodal landscapes\cite{flori2022quantum}. However, given that heuristic methods alone may not guarantee smoothness in the final drug administration policies \cite{wang2020improved}, we refine QPSO’s output using fmincon with a regularization term. This hybrid approach ensures that the resulting treatment strategies not only achieve optimality but also remain clinically feasible by maintaining smooth dosage transitions. Furthermore, we explore different functional forms of drug effectiveness, such as cosine-Gaussian, cosine-sine, and exponential-Gaussian modulation, to assess their impact on tumor suppression. Numerical results indicate that periodic and localized modulation strategies outperform conventional methods, leading to enhanced tumor reduction and reduced drug toxicity. These findings demonstrate the potential of integrating advanced optimization techniques with a deep understanding of tumor dynamics to develop adaptive chemotherapy regimens that improve patient outcomes.

\section{Dynamic Model for Cancer Chemotherapy}

\subsection{Modeling Tumor Heterogeneity and Dynamics}

To model tumor heterogeneity, we introduce a continuous trait $x \in [0,1]$ that represents the normalized level of a cancer cell's drug resistance. The function $n(t,x)$ denotes the population density of cells exhibiting resistance level $x$ at time $t$. Consequently, $n(t,0)$ represents the density of the most drug-sensitive subpopulation, while $n(t,1)$ corresponds to the density of the fully resistant subpopulation. Key biological parameters—including the replication rate $r(x)$, natural death rate $\mu(x)$, and cytotoxic killing parameter $\phi(x)$—are all functions of this trait \cite{norton1977tumor, norton1986norton}. The total tumor cell population, $N(t)$, is the integral of this density across all traits:
\begin{equation}
N(t) = \int_0^1 n(t,x) \, dx.
\end{equation}

For a computationally tractable solution, we discretize the continuous trait space into $m$ distinct traits $x_i$, where $i=1,2,...,m$. This allows us to represent the tumor population as a state vector $\bar{N}(t) \in \mathbb{R}^m$, where each element corresponds to the population of cells with trait $x_i$. The trait-dependent parameters thus become diagonal matrices: the replication rate $R = \text{diag}(r(x_1), \dots, r(x_m))$, the drug toxicity level $\Phi = \text{diag}(\phi(x_1), \dots, \phi(x_m))$, and the natural death rate $M = \text{diag}(\mu(x_1), \dots, \mu(x_m))$. The total tumor population is then approximated by $N(t) = e \bar{N}(t)$, with $e = \tfrac{1}{m}[1,\ldots,1] \in\mathbb{R}^{1\times m}$.

The state equations for trait-specific cancer cells under a single medication are given by \cite{wang2019optimal}:
\begin{equation}
\frac{d}{dt} \bar{N}(t) = \left[ R - \frac{\Phi u}{1 + u} - G(e\bar{N}) M \right] \bar{N},
\end{equation}
and for a double medication, it is:
\begin{align}
\frac{d}{dt} \bar{N}(t) &= \Bigg[ R - \frac{\Phi_1 u_1}{1 + u_1} - \frac{\Phi_2 u_2}{1 + u_2} + \frac{\Phi_{12} u_1}{1 + u_1} \frac{u_2}{1 + u_2} \notag \\
&\quad - G(e \bar{N}) M \Bigg] \bar{N},
\end{align}
Here, $\Phi_1$ and $\Phi_2$ represent the toxicities of the first drug and the second drug, and $\Phi_{12}$ denotes the interaction correction term~\cite{wang2019optimal}. $u$, $u_1$, $u_2$ are drug concentrations, where $0 \leq u, u_1, u_2 \leq u_{\max}$. Here, $G(e \bar{N})=
\log(1+e \bar{N})$ is an increasing function, indicating that as the tumor grows in size, the rate of proliferation of cancer cells decreases while the rate of cell death increases, primarily due to nutrient limitations~\cite{pahnehkolaei2023optimal}.

The objective function $J$, a cost functional to be minimized, is defined for a single drug as:
\begin{equation}
J = \bar{\alpha} \bar{N}(T) + \int_{0}^{T} \left( \bar{\beta} \bar{N}(t) + \gamma u(t) \right) dt,
\end{equation}
and for a double medication as:
\begin{equation}
J = \bar{\alpha} \bar{N}(T) + \int_{0}^{T} \left( \bar{\beta} \bar{N}(t) + \gamma_1 u_1(t) + \gamma_2 u_2(t) \right) dt.
\end{equation}
This functional balances minimizing the tumor burden with the cost of drug administration. The term $\bar{\alpha}\bar{N}(T)$ is the terminal cost for the final tumor size, while the integral accounts for the cumulative tumor burden ($\bar{\beta}\bar{N}(t)$) and drug toxicity ($\gamma u(t)$ or $\gamma_1 u_1(t) + \gamma_2 u_2(t)$) over the treatment period $[0,T]$. The positive weights $\alpha, \beta > 0$ penalize the tumor size, and $\gamma, \gamma_1, \gamma_2 > 0$ penalize drug usage, where $\bar{\alpha} = \alpha e$ and $\bar{\beta} = \beta e$.

For the mutation case, the transition probability $p(x|y)$ for a mutation from trait $y$ to trait $x$ is modeled using a modified Gaussian kernel:
\begin{equation}
p(x|y) = k(y) e^{-\frac{1}{2} \left(\frac{x - y}{\sigma} \right)^2},
\end{equation}
where $\sigma$ is the standard deviation and $k(y)$ is a normalization constant that ensures $\int_0^1 p(x|y)dx=1$ for any given trait $y \in [0,1]$. The replication rate matrix $R$ is then modified by the fraction of cells that mutate, $\theta$ \cite{wang2019optimal2}, and can be expressed as

\begin{equation}
R =
\begin{bmatrix}
r_1(1-\theta) & \theta r_2 p(r_1 | r_2) & \cdots & \theta r_n p(r_1 | r_n) \\
\theta r_1 p(r_2 | r_1) & r_2(1-\theta) & \cdots & \theta r_n p(r_2 | r_n) \\
\vdots & \vdots & \ddots & \vdots \\
\theta r_1 p(r_n | r_1) & \theta r_2 p(r_n | r_2) & \cdots & r_n(1-\theta)
\end{bmatrix}.
\end{equation}
The $R$ matrix models mutation dynamics by adjusting replication rates $r_i$ for each trait $x_i$, with off‑diagonal terms $\theta r_i p(r_j|r_i)$ representing mutation probabilities between traits. The mutation fraction $\theta$ adjusts replication rates based on the likelihood of trait mutates, capturing genetic evolution and simulating tumor adaptation and growth under mutation effects.

\section{Optimal Control Formulation and Hybrid QPSO Solution Method}

In this section, we present a detailed formulation of the optimal control problem for cancer chemotherapy and introduce a hybrid optimization approach that combines QPSO with a gradient-based refinement using \texttt{fmincon}. This hybrid method is designed to address the challenges of high-dimensional optimization in the context of tumor heterogeneity and drug resistance.

\subsection{Step 1: Quantum Particle Swarm Optimization (QPSO)}

QPSO is an advanced variant of the classical PSO algorithm, inspired by quantum mechanics principles. It enhances the exploration capabilities of PSO by modeling particles as quantum entities, allowing them to explore the solution space more efficiently \cite{yang2004quantum}. The pseudocode for QPSO is provided in Algorithm~\ref{alg: QPSO}. The parameters $\phi_{\mathrm{LA}}$, $\beta_{\mathrm{coeff}}$, $r_{\mathrm{a}}$ and $L$ are the local attractor weight, contraction-expansion coefficient, random number uniformly distributed in the interval $(0,1)$, and characteristic length of the potential well, respectively, which guide the exploration of the swarm.

\begin{algorithm}
\caption{Quantum Particle Swarm Optimization (QPSO)}
\label{alg: QPSO}
\begin{algorithmic}[1]
\STATE \textbf{Input:} Population size $N$, max iterations $max\_iter$, bounds $[0, u_{\max}]$, threshold $\epsilon$
\STATE \textbf{Initialize:} Particles $u_i(t)$, personal bests $p_{\text{best}}$, and global best $g_{\text{best}}$.
\FOR{$k=1$ to $max\_iter$}
    \FOR{each particle $i$}
        \STATE Evaluate fitness $J(u_i)$ and update $p_{\text{best}}(i)$ and $g_{\text{best}}$.
        \STATE Compute local attractor: $P_i = \phi_{\mathrm{LA}} p_{\text{best}}(i) + (1 - \phi_{\mathrm{LA}}) g_{\text{best}}$.
        \STATE Update position: $u_i^{(t+1)} = P_i + \beta_{\mathrm{coeff}} L \cdot \text{sign}(r_{\mathrm{a}}-0.5)  \cdot \log(1/(r_{\mathrm{a}}+\epsilon))$.
        \STATE Enforce bounds: $u_i \in [0, u_{\max}]$.
    \ENDFOR
    \STATE Terminate if improvement in $g_{\text{best}}$ is below $\epsilon$.
\ENDFOR
\STATE \textbf{Output:} Optimized dosage schedule $u^*_{\text{QPSO}}(t)$.
\end{algorithmic}
\end{algorithm}

\subsection{Step 2: Regularization using \texttt{fmincon}}

The QPSO algorithm provides a robust global search, but may not guarantee smoothness in control solutions \cite{wang2020improved}. To address this, we employ a gradient-based refinement using \texttt{fmincon}, a nonlinear programming solver. The pseudocode for this step is presented in Algorithm~\ref{alg:fmincon}.

\begin{algorithm}
\caption{Local Refinement using \texttt{fmincon}}
\label{alg:fmincon}
\begin{algorithmic}[1]
\STATE \textbf{Input:} Initial solution from QPSO: $u^*_{\text{QPSO}}(t)$, objective function $J(u)$, constraints $0 \leq u \leq u_{\max}$
\STATE \textbf{Define optimization problem:}
\STATE \quad Minimize $J(u)$ subject to $0 \leq u \leq u_{\max}$
\STATE \textbf{Call \texttt{fmincon} solver}: Solve for refined dosage schedule:
$u^*_{\text{final}} = \arg \min\limits_{u} J(u)$
\STATE \quad Ensure smoothness and adherence to clinical constraints
\STATE \textbf{Output:} Optimized dosage schedule $u^*_{\text{final}}(t)$
\end{algorithmic}
\end{algorithm}

This hybrid approach first employs QPSO for a global search, identifying promising drug dosage schedules. Then \texttt{ fmincon} refines these solutions taking $u^*_{\text{QPSO}}(t)$ as initial values, ensuring smooth and clinically viable treatment regimens. Using global and local optimization techniques, this method improves the robustness and effectiveness of cancer chemotherapy planning.

\section{Verification: Optimality Condition}

To verify the optimality of the obtained control solutions, we perform first- and second-order optimality checks using the Hamiltonian function. For the single-medication case, the Hamiltonian is formulated as:

\begin{equation}
H = \bar \beta
\bar N(t) + \gamma u(t) + \lambda^T \left[ R - \frac{\Phi u}{1 + u} - G(e \bar{N}) M \right] \bar{N},
\end{equation}
where $\bar\beta=\beta e$, with $ \beta=400$ and $\gamma=4000$ are weighting coefficients, $\lambda$ is the costate vector, $R$ is the replication rate, $\Phi$ is the cytotoxic killing parameter, $G(e \bar{N})$ is the growth function, and $M$ is the natural death rate.

\begin{figure}[!t]
    \centering
    \includegraphics[width=\linewidth, height=6cm]{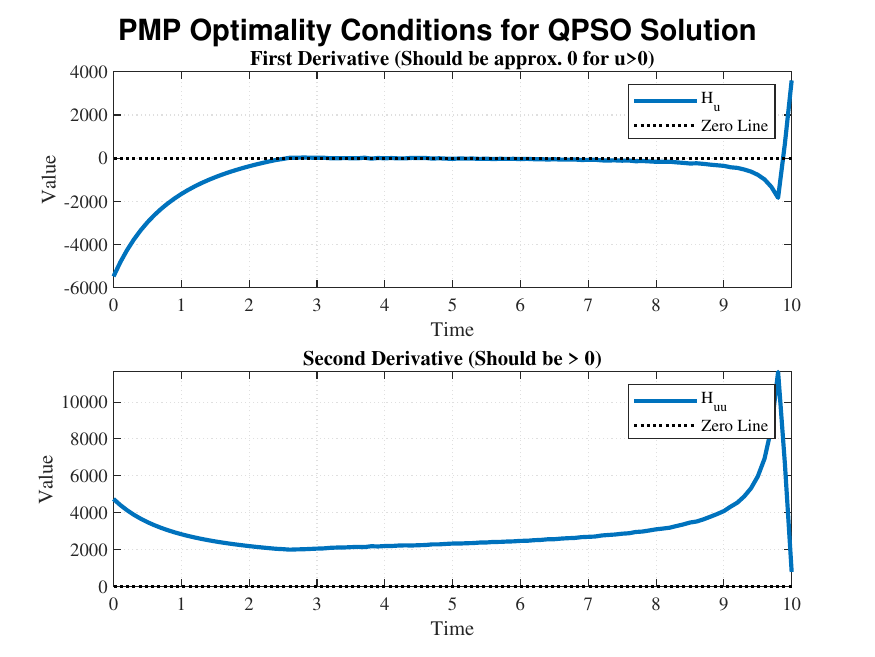}
    \caption{Optimality check results: (a) First-order check showing the gradient of the Hamiltonian with respect to the control variable, and (b) Second-order check confirming the positive definiteness of the second derivative, ensuring Hamiltonian minimization.}
    \label{fig:optimality_check}
\end{figure}

The costate equation, which describes the dynamics of the costate vector, is derived from $\dot{\lambda} = - \frac{\partial H}{\partial \bar{N}}$. The first-order optimality check (necessary condition) ensures that the gradient of the Hamiltonian with respect to the control variable is zero:

\begin{equation}
\frac{\partial H}{\partial u} = \gamma - \frac{\lambda^T \Phi \bar{N}}{(1 + u)^2} = 0.
\end{equation} The second-order optimality check (sufficient condition) ensures that the Hamiltonian is minimized by verifying that the second derivative is positive definite:

\begin{equation}
\frac{\partial^2 H}{\partial u^2} = \frac{2 \lambda^T \Phi \bar{N}}{(1 + u)^3}.
\end{equation} Figure \ref{fig:optimality_check} demonstrates the optimality of the control process for a single medication, non-mutation system where $R=\text{diag}\left[\frac{2}{1+3x^4}\right]$, $\Phi=\text{diag}\left[-\sin(x-1)+1.5\right]$, and $M = 0.5\,I_n$\cite{pahnehkolaei2023optimal}. Optimality is confirmed through first-order and second-order conditions, subject to the constraint $u \in [0, u_{max}=3]$. In Figure 1(a), deviations of $\frac{\partial H}{\partial u}$ from zero in certain regions are attributed to the presence of constraints, which affect the control trajectory and enforce feasibility. Figure 1(b) shows that the second-order derivative of the Hamiltonian is positive definite, ensuring the optimality of the solution.

\section{Numerical Examples}
In this section, we investigate the effectiveness of double-medication strategies for controlling tumor populations, considering different functional forms for drug effectiveness. The goal is to minimize tumor burden while ensuring a uniform distribution across different tumor traits. We analyze three distinct combination therapy scenarios to determine how different modulation strategies impact treatment efficacy.

\subsection{Model Parameters}

The final tumor population weight is $\alpha = 5$, and the intermediate weight is $\beta = 400$. Toxicity weights for Drug 1 and Drug 2 are $\gamma_1 = 4000$ and $\gamma_2 = 4000$, respectively. The death rate matrix is $M = 0.5\,I_n$, with 21 equally spaced trait values between 0 and 1, providing sufficient resolution to capture nonlinear variations while maintaining computational efficiency. The replication rate is $r = \frac{2}{1 + 3x}$, and the maximum drug concentrations for Drug 1 and Drug 2 are $u_{\text{max1}} = 3$ and $u_{\text{max2}} = 3$. The interaction coefficient is $\phi_{12} = 0.2$, with mutation parameters $\sigma = 5$ and $\theta = 0.2$. The time variable spans a duration of 10 years, during which the tumor population dynamics and drug dosages are simulated.

Building on the use of qualitative functions for the cytotoxic killing parameter, $\phi(x)$ \cite{wang2016optimal}, we move beyond prior studies that relied mainly on polynomial or rational forms \cite{greene2014impact, schattler2015optimal}. We introduce localized functions (e.g., Gaussian) to model targeted therapies against narrow phenotypes \cite{dagogo2018tumour}, and periodic functions (e.g., cosine) to capture agents with non-monotonic effects across resistance. This broader class of functions allows a more complete assessment of how drug-effectiveness profiles shape optimal control.

\subsection{Combination Therapy Scenarios}
We analyze three distinct combination therapy strategies, each characterized by a unique modulation of drug effectiveness. These strategies determine how cytotoxic effects are distributed across tumor traits, influencing both tumor suppression and the final tumor distribution. The results for each case include the optimal drug administration profiles, tumor population evolution, and the final tumor trait distribution.

\subsubsection{Cosine-Gaussian Combination}
In this scenario, Drug 1 follows a cosine-squared function, providing broad and periodic suppression across tumor traits, while Drug 2 follows a Gaussian distribution, concentrating its effectiveness on a specific subset of trait values. The respective functional forms are:

\begin{equation}
    \left\{
    \begin{aligned}
        \phi_1(x) &= 1 + \cos^2(0.25\pi x) \\
        \phi_2(x) &= 1 + \exp\left(-\frac{(x-0.2)^2}{2 \times 0.25^2}\right)
    \end{aligned}
    \right.
\end{equation}

\begin{figure}[!t]
    \centering
    \includegraphics[width=\linewidth, height=6cm]{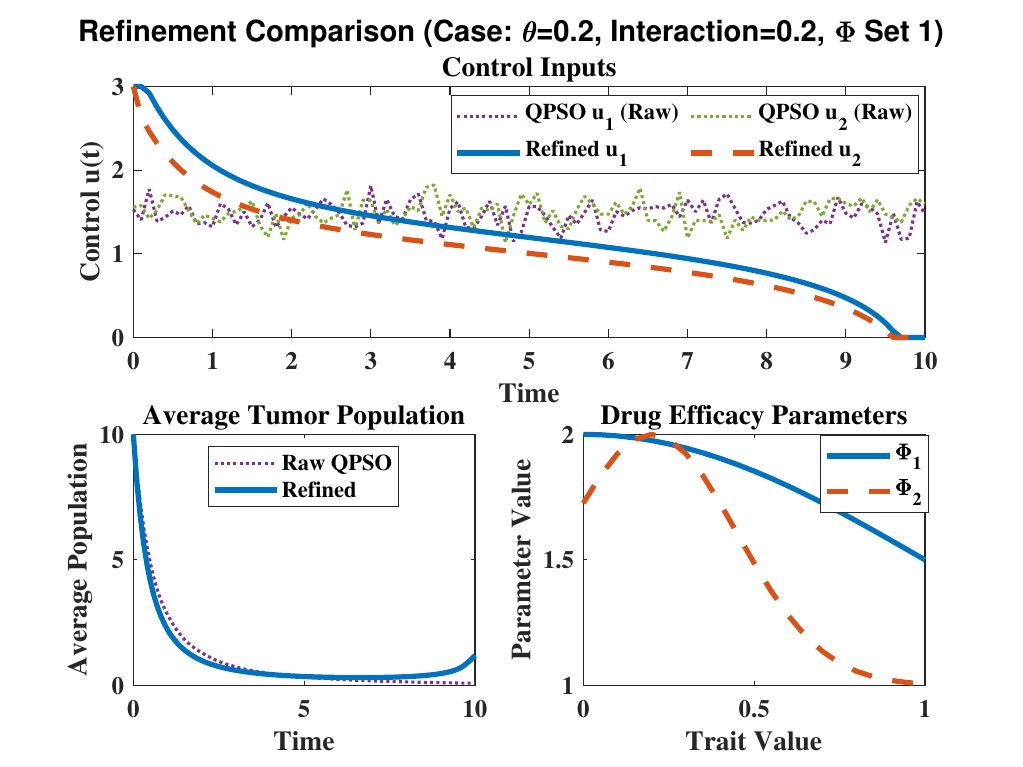}
    \caption{Benefit of the two-step hybrid optimization. (Top) Comparison of the raw QPSO control (dashed) and the final refined control (solid). (Bottom Left) Resulting tumor population dynamics. (Bottom Right) Cytotoxic killing parameters.}
    \label{fig:cosine_gaussian}
\end{figure}

The optimal control strategy, shown in Figure~\ref{fig:cosine_gaussian}, indicates that Drug 1 is initially administered at a high dosage, gradually decreasing over time. Drug 2 follows a similar but lower trajectory, suggesting a complementary interaction. Drug 1 provides broad suppression across tumor traits, while Drug 2 selectively targets subpopulations more susceptible to its cytotoxic effects. Figure~\ref{fig:cosine_gaussian} also illustrates the tumor population dynamics under this regimen. A sharp initial decline in tumor burden is observed, reflecting effective early suppression. However, a minor resurgence occurs in later stages, likely due to the emergence of resistant subpopulations. The cytotoxic killing response reveals that Drug 1 maintains moderate effectiveness across trait values, whereas Drug 2 exhibits a pronounced peak in cytotoxicity at intermediate trait values, reinforcing selective eradication.

To demonstrate the benefit of our hybrid method, Fig. \ref{fig:cosine_gaussian} compares the raw control trajectory from QPSO with the final refined trajectory. The raw output (dashed lines) identifies the optimal control landscape but has high-frequency oscillations unsuitable for clinical use. The regularization step smooths these into a practical schedule (solid lines). Importantly, this refinement is achieved with a negligible impact on therapeutic efficacy, as shown by the nearly identical tumor population dynamics.

\subsubsection{Cosine-Sine Contrast}
This case explores an alternating suppression strategy where Drug 1 follows a cosine-squared function, while Drug 2 adopts a sine-squared function. This ensures that different tumor traits are targeted at different phases of treatment:

\begin{equation}
    \left\{
    \begin{aligned}
        \phi_1(x) &= 1 + \cos^2\left(\frac{\pi}{4} x\right) \\
        \phi_2(x) &= 1 + \sin^2\left(\frac{\pi}{4} x\right)
    \end{aligned}
    \right.
\end{equation}

\begin{figure}[!t]
    \centering
    \includegraphics[width=\linewidth, height=6cm]{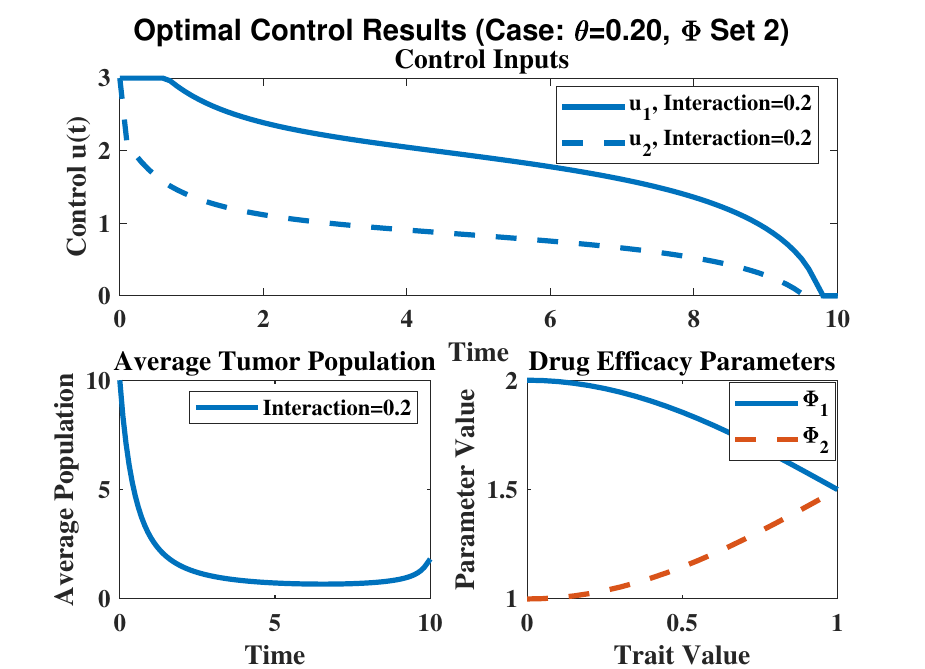}
    \caption{Tumor suppression dynamics under Cosine–Sine therapy. Top: optimal control input; bottom left: tumor population evolution; bottom right: cytotoxic killing parameter. Alternating modulation shifts drug pressure across traits—Drug 1 (cosine-squared) targets lower traits, while Drug 2 (sine-squared) targets higher traits—yielding balanced suppression.}
    \label{fig:cosine_sine}
\end{figure}

Figure~\ref{fig:cosine_sine} shows that Drug 1 begins at a high dose and tapers, while Drug 2 starts lower but remains sustained, reinforcing suppression later. Tumor dynamics display a sharp initial decline followed by stabilization, with a modest late-stage rise indicating resistant subpopulations. The killing response transitions smoothly between drugs, with Drug 1 dominating at lower traits and Drug 2 compensating at higher traits.

\subsubsection{Exponential-Gaussian Mix}
In this scenario, Drug 1 follows an exponential decay function, targeting lower trait values early in treatment, while Drug 2 follows a Gaussian function, concentrating its cytotoxic effects on intermediate traits:

\begin{equation}
    \left\{
    \begin{aligned}
        \phi_1(x) &= 1 + \exp(-3x) \\
        \phi_2(x) &= 1 + \exp\left(-\frac{(x-0.5)^2}{2 \times 0.25^2}\right)
    \end{aligned}
    \right.
\end{equation}

\begin{figure}[!t]
    \centering
    \includegraphics[width=\linewidth, height=6cm]{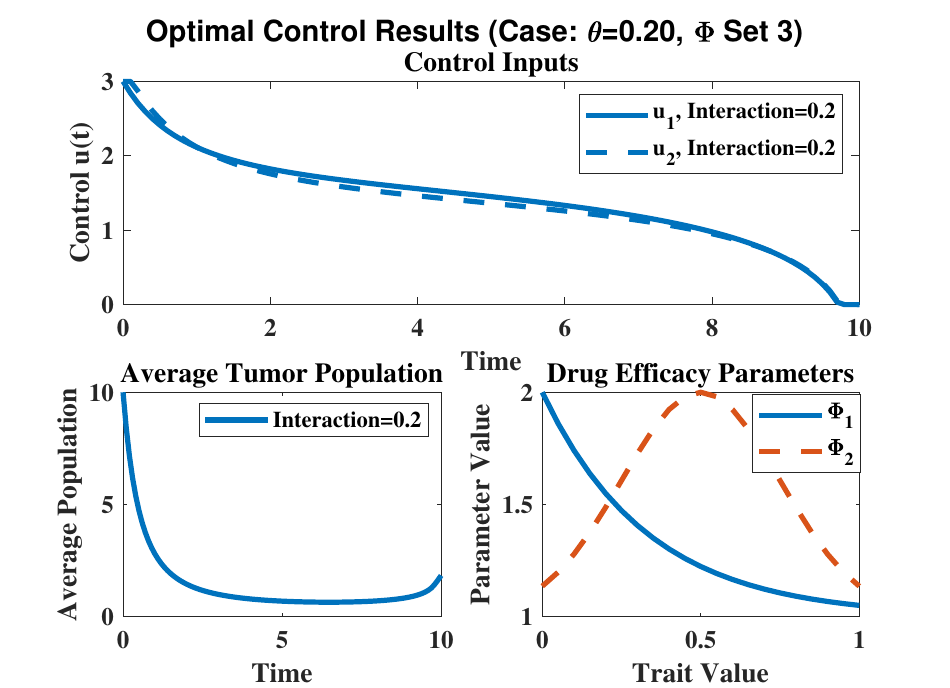}
    \caption{Tumor suppression dynamics under Exponential–Gaussian therapy. Top: optimal control input; bottom left: tumor population evolution; bottom right: cytotoxic killing parameters. Drug 1 drives early suppression, followed by selective action from Drug 2, producing localized resistance at lower trait values.}
    \label{fig:exponential_gaussian}
\end{figure}

As shown in Figure~\ref{fig:exponential_gaussian}, Drug 1 is delivered at a high initial dose that tapers over time, while Drug 2 follows a slower but sustained course. Their combination achieves broad suppression with targeted reinforcement at specific traits. Tumor dynamics reveal a sharp early decline followed by late-stage resurgence in lower-trait subpopulations, suggesting resistance. The killing response reflects a Gaussian-like modulation for Drug 2, whereas Drug 1’s effect decreases monotonically, supporting selective eradication.

\begin{figure}[!t]
    \centering
    \includegraphics[width=\linewidth, height=6cm]{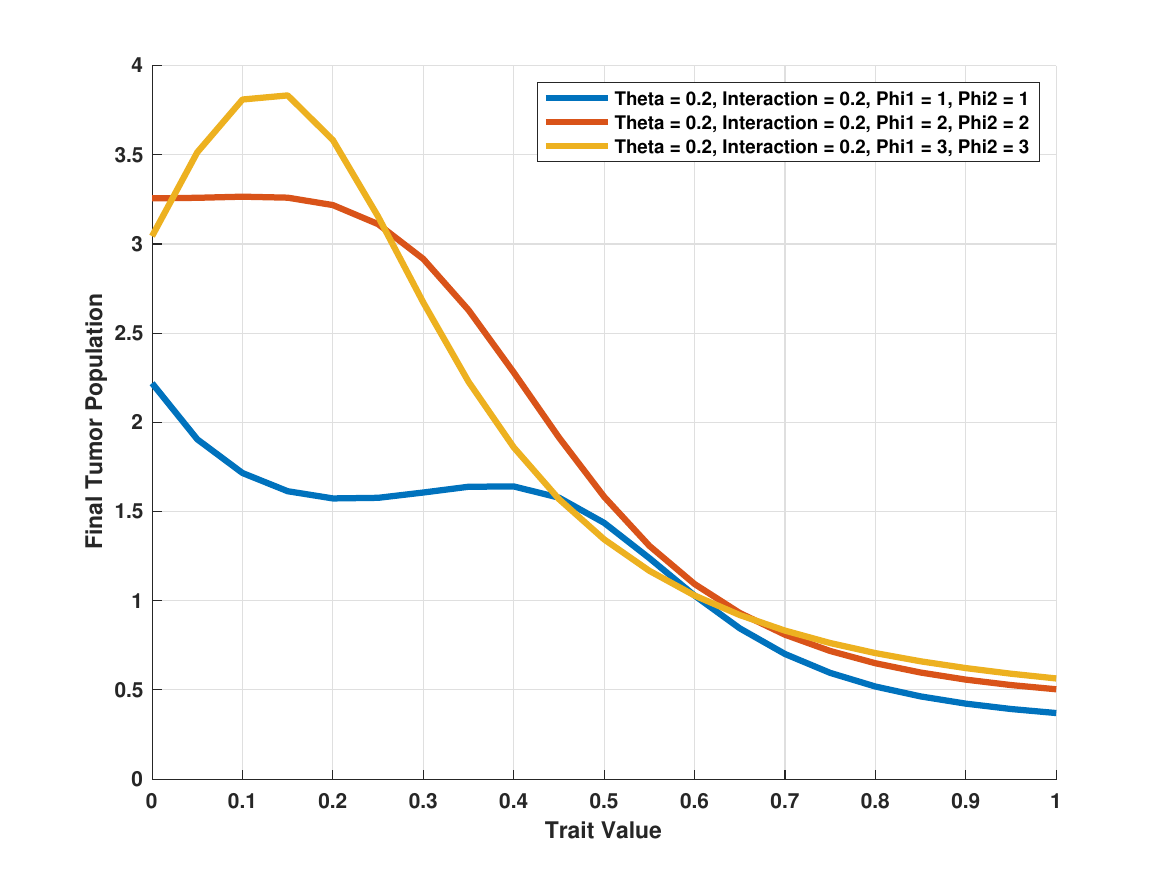}
    \caption{Final tumor distributions for all therapy strategies. (a) Cosine-Gaussian (blue), (b) Cosine-Sine (red), and (c) Exponential-Gaussian (yellow). The Cosine-Gaussian therapy achieves uniform suppression, while the Cosine-Sine approach results in intermediate trait resistance. The Exponential-Gaussian strategy shows localized resistance at lower trait values, indicating potential treatment limitations.}
    \label{fig:final_comparison}
\end{figure}

The final tumor distributions for all therapy strategies are compared in Figure~\ref{fig:final_comparison}. The Cosine-Gaussian approach results in a uniform decrease in trait values, while the Cosine-Sine therapy exhibits intermediate trait resistance. The Exponential-Gaussian strategy shows localized resistance at lower trait values, highlighting the potential limitations of this approach in eradicating all tumor subpopulations.

\begin{figure}[!t]
    \centering
    \includegraphics[width=\linewidth, height=6cm]{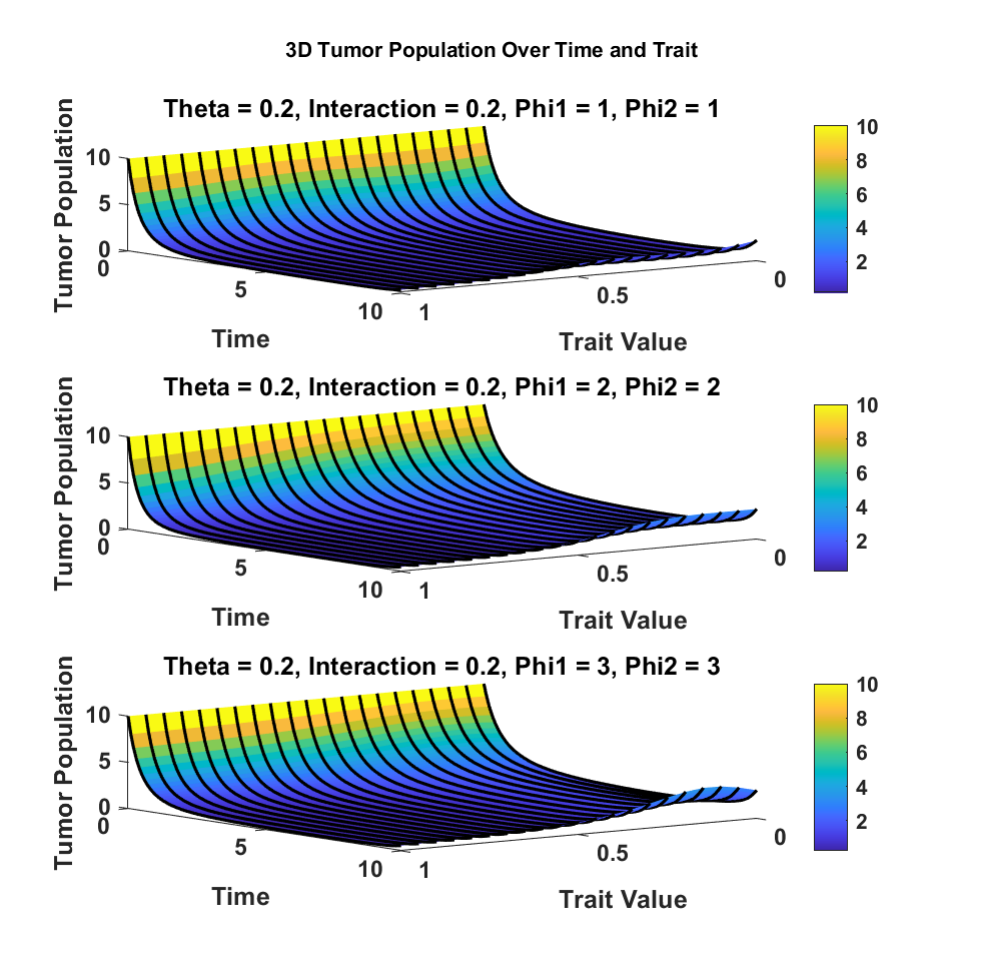}
    \caption{Tumor population evolution under different treatment regimens. (top) Cosine-Gaussian, (middle) Cosine-Sine, and (bottom) Exponential-Gaussian.}
    \label{fig:3d_plot}
\end{figure}

 Figure~\ref{fig:3d_plot} presents the temporal evolution of the tumor population under each therapy regimen. Periodic suppression in Cosine-Gaussian and Cosine-Sine therapies maintains consistent tumor reduction. In contrast, Exponential Gaussian therapy exhibits transient resistance at lower trait values, emphasizing the importance of selecting modulation strategies that effectively counter tumor heterogeneity.

\section{Discussion}
The functional form of drug effectiveness significantly influences tumor suppression and the final distribution of traits. Periodic functions, such as the cosine-based functions examined here, distribute drug impact more evenly across tumor traits, proving effective for heterogeneous tumor populations. This is evident in the Cosine-Gaussian and Cosine-Sine scenarios, where the final tumor distributions are smooth (Figure~\ref{fig:final_comparison}), indicating consistent tumor control. In contrast, the Exponential-Gaussian combination results in a non-uniform impact, with a pronounced peak in the final tumor population. This suggests that while the exponential-Gaussian approach targets specific traits effectively, it may allow other subpopulations to persist.

Combination therapy strategies incorporating periodic cytotoxic functions enhance tumor suppression beyond single-drug treatments and offer greater control over tumor heterogeneity. They also reduce overall toxicity by maintaining efficacy while limiting excessive drug exposure. However, the mathematical model used does not fully capture tumor biology and drug interactions, necessitating careful parameter tuning and experimental validation for clinical application. Future research should focus on refining drug interaction models, exploring adaptive control strategies, explicitly modeling drug resistance evolution, and leveraging advanced numerical optimization techniques for robust and personalized drug administration schedules.
A key motivation for this work was to develop a method robust to initial conditions, a known issue with PMP-based approaches. Our QPSO framework does not rely on costate integration and instead uses a population-based global search. This structure ensures reliable convergence to a high-quality solution without the fragility associated with selecting a precise initial costate guess.

\section{CONCLUSION}

This study introduces a hybrid QPSO-based optimal control framework for cancer chemotherapy that tackles tumor heterogeneity and mutation dynamics. By combining QPSO’s global search with regularization-based refinement, it overcomes the sensitivity and convergence issues of traditional costate methods. Simulations highlight the effectiveness of periodic dosing, especially cosine-based schedules, in limiting resistance and ensuring uniform suppression. The results emphasize the role of spatial drug distribution, while future work should integrate patient-specific traits, complex drug interactions, and refined resistance models to enhance clinical relevance.

\end{document}